\documentclass[12pt]{article}
\usepackage{amsxtra,amssymb,amsthm,amsmath,latexsym}

\theoremstyle{plain}

\def\R{{\mathbb R}}

\def\oH{\buildrel\circ\over H}
\def\oH1{\buildrel\circ\over H\kern-.02in{}^1}

\begin{document}

\title{ A symmetry problem
   \thanks{key words: symmetry for PDE, harmonic functions 
    }
}

\author{
A.G. Ramm\\
Mathematics Department, 
Kansas State University, \\
 Manhattan, KS 66506-2602, USA\\
ramm@math.ksu.edu\\
}

\date{}

\maketitle\thispagestyle{empty}

\begin{abstract}
The following result is proved:

{\bf Theorem.} Let $D\subset \R^3$ be a bounded domain
homeomorphic to a ball, $|D|$ be its volume, $|S|$ be the
surface area of its smooth boundary $S$, $D\subset
B_R:=\{x:|x|\leq R\}$, and $H_R$ is the set of all harmonic
in $B_R$ functions. If $$\frac 1 {|D|}\int_Dhdx=\frac 1
{|S|}\int_Shds\quad \forall h\in H_R,$$ then $D$ is a ball.  
\end{abstract} 
\section{Introduction}
Consider a bounded domain $D \subset \R^n$, $n = 3,$ with a
boundary $S$. The exterior domain is $D^\prime = \R^3
\backslash D$. Assume that $S$ is Lipschitz. Let $N$ be the
exterior unit normal to $S$. Our result is:

{\bf Theorem.} {\it Let $D$ be a bounded domain homeomorphic 
to a ball,
$|D|$ be its volume, $|S|$ be the surface area of 
its smooth boundary
$S$,  $D\subset B_R:=\{x:|x|\leq R\}$, 
and $H_R$ is the set of all
harmonic in $B_R$ functions. If
$$\frac 1 {|D|}\int_Dhdx=\frac 1 {|S|}\int_Shds\quad \forall h\in H_R,
\eqno{(1)}$$
then $D$ is a ball.}

We prove this theorem using the result from [1] and the 
following lemmas:

{\bf Lemma 1.} {\it The set $\{h|_S\}_{\forall h\in H_R}$ is dense in 
$L^2(S)$.} 

For convenience of the reader let us formulate the result from [1]:

{\bf Lemma 2.} {\it If $D$ is homeomorphic to a ball and  
$$\Delta u= \frac 1 {|D|} \text { in } D,\quad u|_S=0, \quad 
u_N|_S=const,$$
then $D$ is a ball.} 

The result is valid for $D\in \R^n,\, n\geq 2$, and its 
proof is essentially the same. If $D$ is a ball, then 
formula (1) holds and both integrals are equal to $h(0)$,
where the origin is at the center of the ball (mean-value 
theorem for harmonic functions).

In the next Section proofs are given.

\section{Proofs} 

{\it Proof of Lemma 1.} Assume the contrary. Then there exists an $f$ such
that $\int_S fhds=0$ for all $h\in H_R$. Take $h=\frac 1 {4\pi |x-y|}$
where $y\in B_R':=\R^3\setminus B_R$ is arbitrary. Then
$$w(y):=\int_S\frac {fds}{4\pi |s-y|}=0\quad \forall y\in B_R'.$$
Since $w$ is harmonic in $D'$, it follows that $w=0$ in $D'$.
Since $w$ is harmonic in $D$, it follows that $w=0$ in $D$.
By the jump formula for the normal derivative of the simple-layer 
potential, one concludes that $f=0$.
Lemma 1 is proved. $\Box$.

{\it Proof of Theorem 1.} Let $u$ be the unique solution to 
the problem:
$$\Delta u=\frac 1 {|D|} \quad \text { in } D, \quad u|_S=0.
 \eqno{(2)}$$

 Multiply (2) by an arbitrary $h\in H_R$, integrate over $D$, then by 
parts, use the boundary condition (2) and get:
$$\frac 1 {|D|}\int_Dhdx=\int_Su_Nhds\quad \forall h\in H_R.
\eqno{(3)}$$
By (1) this implies:
$$\int_S(u_N-\frac 1 {|S|})hds=0\quad \forall h\in H_R.
\eqno{(4)}$$
By Lemma 1, equation (4) implies 
$$
u_N=\frac 1 {|S|} \quad \text {on } S.
\eqno{(5)}$$
By Lemma 2, equations (2) and (5) imply that $D$ is a ball.
Theorem 1 is proved. $\Box$

\end{document}